\documentclass[11pt,letterpaper]{amsart}
\usepackage{mathrsfs,graphicx}
\usepackage{amsopn,amssymb,bbold}
\usepackage[latin1]{inputenc}

\newcommand{\Z}{\mathbb{Z}}
\newcommand{\R}{\mathbb{R}}
\newcommand{\C}{\mathbb{C}}

\newcommand{\E}{\mathscr{E}}

\newcommand{\rhoth}{\rho_{\text{Th}}}

\newcommand{\T}{\mathscr{T}}
\renewcommand{\P}{\mathscr{P}}

\newcommand{\CP}{\mathbb{CP}}
\newcommand{\F}{\mathscr{F}}

\renewcommand{\H}{\mathbb{H}}

\newcommand{\ML}{\mathscr{M\!L}}
\newcommand{\PML}{\mathbb{P} \! \mathscr{M\!L}}
\newcommand{\MF}{\mathscr{M\!F}}
\newcommand{\PMF}{\mathbb{P} \! \mathscr{M\!F}}
\newcommand{\PQ}{\mathbb{P}^+Q}

\newcommand{\comment}[1]{}

\DeclareMathOperator{\aut}{Aut}
\DeclareMathOperator{\PSL}{\mathrm{PSL}}

\DeclareMathOperator{\gr}{gr}
\DeclareMathOperator{\pr}{pr}
\DeclareMathOperator{\cl}{\kappa}
\DeclareMathOperator{\Gr}{Gr}

\newcommand{\boldpoint}[1]{\medskip\par\noindent\textbf{#1}}
\newcommand{\sboldpoint}[1]{\smallskip\par\noindent\textbf{#1}}

\DeclareMathOperator{\im}{Im}
\DeclareMathOperator{\re}{Re}

\theoremstyle{plain}

\newtheorem{thm}{Theorem}[section]
\newtheorem{cor}[thm]{Corollary}
\newtheorem{conj}[thm]{Conjecture}
\newtheorem{lem}[thm]{Lemma}
\newtheorem{prop}[thm]{Proposition}
\theoremstyle{definition}

\theoremstyle{remark}
\newtheorem*{rem}{Remark}

\newtheorem*{note}{Note}
\newtheorem*{notes_env}{Notes}

\theoremstyle{plain}
\newtheorem*{blankthm}{Theorem}

\begin{document}

\title[Grafting, Pruning, and the Antipodal Map]{Grafting, Pruning,
and the Antipodal Map on Measured Laminations}
\author{David Dumas}
\date{November 30, 2006}
\address{Department of Mathematics \\ Brown University \\ Providence, RI 02912 \\ USA\\}
\email{ddumas@math.brown.edu}
\thanks{The author was partially supported by an NSF Postdoctoral Research Fellowship}
\begin{abstract}
Grafting a measured lamination on a hyperbolic surface defines a
self-map of Teichmüller space, which is a homeomorphism by a result
of Scannell and Wolf.  In this paper we study the large-scale behavior
of pruning, which is the inverse of grafting.

Specifically, for each conformal structure $X \in \T(S)$, pruning $X$
gives a map $\ML(S) \to \T(S)$.  We show that this map extends
to the Thurston compactification of $\T(S)$, and that its boundary
values are the natural \emph{antipodal involution} relative to $X$ on
the space of projective measured laminations.

We use this result to study Thurston's grafting coordinates on the
space of $\CP^1$ structures on $S$.  For each $X \in \T(S)$,
we show that the boundary of the space $P(X)$ of $\CP^1$ structures on
$X$ in the compactification of the grafting coordinates is the graph
$\Gamma(i_X)$ of the antipodal involution $i_X : \PML(S) \to \PML(S)$.
\end{abstract}

\maketitle

\tableofcontents

\section{Introduction}

Grafting is a procedure that begins with a hyperbolic structure $Y \in
\T(S)$ in the Teichmüller space of a surface $S$ of negative Euler
characteristic and a measured geodesic lamination $\lambda \in
\ML(S)$.  By replacing $\lambda$ with a thickened version that carries
a natural Euclidean metric, a new conformal structure $X = \gr_\lambda
Y \in \T(S)$, the \emph{grafting of $Y$ along $\lambda$}, is obtained.

Scannell and Wolf have shown that for
each lamination $\lambda \in \ML(S)$, the conformal grafting map
$\gr_\lambda : \T(S) \to \T(S)$ is a homeomorphism, thus there is an
inverse or \emph{pruning map} $\pr_\lambda : \T(S) \to \T(S)$
\cite{SW}.

In this paper we describe the large-scale behavior of pruning $X \in
\T(S)$ in terms of the conformal geometry of $X$.  This description is
based on the map
$$ \Lambda : Q(X) \to \ML(S) $$
which records the measured lamination equivalent to the horizontal
foliation of a holomorphic quadratic differential.  Hubbard and Masur
showed that $\Lambda$ is a homeomorphism \cite{HM}, so we can use it to transport the
involution $(\phi \mapsto -\phi)$ of $Q(X)$ to an involutive homeomorphism $i_X :
\ML(S) \to \ML(S)$.  Since the map $\Lambda$ is homogeneous, $i_X$ descends
to an involution on $\PML(S) = (\ML(S) - \{0\}) / \R^+$,
$$ i_X : \PML(S) \to \PML(S),$$ 
which we call the \emph{antipodal involution with respect to $X$},
since it is conjugate by $\Lambda$ to the actual antipode map on the
vector space $Q(X)$.

\boldpoint{Statement of results.}  In §\ref{sec:mainproof}
we show that $i_X$ governs the large-scale behavior of the map $\ML(S)
\to \T(S)$ which prunes $X$ along a given lamination, i.e.
$\lambda \mapsto \pr_\lambda X$.  Specifically, let
$\overline{\ML(S)}$ denote the natural compactification of $\ML(S)$ by
$\PML(S)$, and let $\overline{\T(S)}$ denote the Thurston
compactification of Teichmüller space, which also has boundary
$\PML(S)$.

\begin{thm}[Antipodal limit]
\label{thm:limit}
The pruning map with basepoint $X$, written $\lambda \mapsto \pr_\lambda X$,
extends continuously to a map $\overline{\ML(S)} \to \overline{\T(S)}$
whose boundary values are exactly the antipodal map $i_X : \PML(S) \to
\PML(S)$.
\end{thm}

An equivalent formulation of Theorem \ref{thm:limit} characterizes the
fibers of the grafting map $\gr: \ML(S) × \T(S) \to \T(S)$:

\begin{thm}[Fibers of grafting]
\label{thm:fiber}
Let $M_X = \gr^{-1}(X) \subset \ML(S) × \T(S)$ denote the set of
all pairs $(\lambda,Y)$ that graft to give $X \in \T(S)$.  Then the
boundary of $M_X$ in $\overline{\ML(S)} × \overline{\T(S)}$ is
the graph of the antipodal involution, i.e.
$$\overline{M_X} = M_X \sqcup \Gamma(i_X)$$ where
$$\Gamma(i_X) = \{ ([\lambda],[i_X(\lambda)]) \; | \; \lambda \in
\ML(S) \} \subset \PML(S) × \PML(S)$$
\end{thm}

The equivalence of Theorems \ref{thm:limit} and \ref{thm:fiber}
follows from the definition of pruning (see §\ref{sec:grafting}).

Let $\P(S)$ denote the space of $\CP^1$ surfaces, i.e. surfaces with
geometric structures modeled on $(\CP^1,\PSL_2(\C))$, and $P(X) \subset
\P(S)$ the subset having underlying conformal structure $X$.
Thurston's parameterization of $\CP^1$ structures is based on a homeomorphism
$$ \ML(S) × T(S) \xrightarrow{\sim} \P(S) $$ that is a generalization of
grafting (see \cite{KT}).  This allows us to view $P(X)$ as a subset of $\ML(S) ×
\T(S)$, and it is natural to try to understand the relationship
between the conformal structure $X$ and the pairs $(\lambda,Y) \subset
\ML(S)×\T(S)$ thus determined.  Equivalently, it is natural to compare
Thurston's geometric parameterization of $\CP^1$ surfaces via grafting
with the classical complex-analytic viewpoint.

Theorem \ref{thm:fiber} gives some information about how
these two perspectives on $\P(S)$ compare, as it provides an explicit 
description of the boundary of $\overline{P(X)}$ as a subset
of $\overline{\ML(S)} × \overline{\T(S)}$:

\begin{thm}[Boundary $P(X)$]
\label{thm:boundary}
For each $X \in \T(S)$, the boundary of $P(X)$ is the graph of the
antipodal involution $i_X$:
$$ \partial P(X) = \Gamma(i_X) \subset \partial (\overline{\ML(S)}
× \overline{\T(S)})$$
\end{thm}

In particular, the closure of $P(X)$ is a ball of dimension $6g-6$,
where $g$ is the genus of $S$.  Theorem \ref{thm:boundary} follows
immediately from Theorem \ref{thm:fiber} because $P(X)$ corresponds to
a fiber Thurston's $\CP^1$-structure generalization of grafting (see
§\ref{sec:grafting}).

\boldpoint{Collapsing and harmonic maps.}  Our main results are obtained by studying the
\emph{collapsing} map $\cl : \gr_\lambda Y \to Y$ that
collapses the grafted portion of the surface back to the geodesic
lamination on $Y$.

The key observation that drives the proof of Theorem \ref{thm:limit}
is that the collapsing map is nearly harmonic, having energy exceeding
that of the associated harmonic map by at most a constant depending
only on the topology of $S$.  This energy comparison is due to
Tanigawa \cite{Ta}, and relies on an inequality of Minsky for length
distortion of harmonic maps \cite{Mi}.

\boldpoint{Geometric limits.}  For the purposes of Theorem
\ref{thm:limit}, which is a result about the asymptotic behavior of
pruning, we are interested in pruning a fixed surface along a
divergent sequence of laminations.  For such a sequence, the energies
of the associated collapsing maps also diverge, while remaining
withing $O(1)$ of the energies of the homotopic harmonic maps.  This
allows us to promote the collapsing map to a genuine harmonic map by
rescaling the hyperbolic metric of $Y$ and taking a geometric limit.
This limit construction relies on the compatibility between
equivariant Gromov-Hausdorff convergence and harmonic maps from
Riemann surfaces as established by Bestvina \cite{Be} and Paulin
\cite{Pa}; the more general theory of Korevaar and Schoen (see
\cite{KS}, \cite{KS2}) for harmonic maps to metric spaces could also
be used.

The resulting limit harmonic map has values in an $\R$-tree, and there
is a detailed structure theory for such harmonic maps due to Wolf
(\cite{W}, \cite{W5}, \cite{W6}).  It follows that a single
holomorphic quadratic differential--the Hopf differential of the limit
harmonic map--encodes both the Thurston limit of the pruned surfaces
and the projective limit of the grafting laminations via its vertical
and horizontal foliations.  This gives rise to the antipodal
relationship expressed by Theorem \ref{thm:limit}.

\boldpoint{Asymmetry of Teichmüller geodesics.}  
It is tempting to compare the role of the antipodal involution in
Theorem \ref{thm:limit} to that of the geodesic involution in a
symmetric space.  Indeed, the antipodal map relative to $X$ exchanges
projective measured laminations $([\lambda],[\mu])$ defining
Teichmüller geodesics that pass through $X$, just as the geodesic
involution exchanges endpoints (at infinity) of geodesics through a
point in a symmetric space.

However, in an appendix we provide an example showing that this
analogy does not work, because Teichmüller geodesics are badly behaved
with respect to the Thurston compactification.  Our example is based
on properties of Teichmüller geodesic rays associated to Strebel
differentials established by Masur (see \cite{Masur2}, \cite{Masur}).
Specifically, we construct a pair of Teichmüller geodesics through a
single point in $\T(S)$ that are asymptotic to each other in one
direction, but which have distinct limit points in the other
direction.  This precludes the existence of any map of the Thurston
boundary that plays the role of a Teichmüller geodesic involution.

\boldpoint{Outline of the paper.}
\sboldpoint{Section \ref{sec:grafting}} introduces the main objects we study in
this paper--the grafting, pruning, and collapsing maps.  Using the
definitions and basic properties of these objects, the equivalence of
Theorems \ref{thm:limit}-\ref{thm:boundary} is explained.

\sboldpoint{Sections \ref{sec:tensors}-\ref{sec:rtrees}} contain further background
material needed for the proof of the main theorem.  Specifically,
Section \ref{sec:tensors} contains definitions and notation related to
quadratic differentials (holomorphic and otherwise) and conformal
metrics that are used in comparing the collapsing map to a harmonic
map.  The antipodal map on measured laminations is then introduced in
Section \ref{sec:foliations} by means of the measured foliations
attached to holomorphic quadratic differentials.  Section
\ref{sec:rtrees} then describes the $\R$-trees dual to measured
foliations and measured laminations that appear naturally when
studying limits of harmonic maps.

\sboldpoint{Section \ref{sec:harm}} discusses results of Wolf (from \cite{W4},
\cite{W5}, \cite{W6}) on harmonic maps between surfaces and $\R$-trees
that allow us to characterize the extension of pruning in terms of the
antipodal map on measured laminations.

\sboldpoint{Section \ref{sec:energy}} presents the key energy
estimates that allow us to compare the collapsing map (and a dual
object, the co-collapsing map) to an associated harmonic map.  The
energy of the collapsing map was first computed by Tanigawa \cite{Ta}
in order to prove that grafting is proper; here we provide an
analogous result on the properness of pruning.

\sboldpoint{Section \ref{sec:hopf}} presents a complementary study of
the Hopf differentials of the collapsing an co-collapsing maps.  Here
we establish a relationship between these (non-holomorphic) quadratic
differentials and the grafting lamination that is analogous to the
relationship between a holomorphic quadratic differential and its
associated measured foliation.  This relationship is essential in the
proof of the main theorem.

\sboldpoint{Section \ref{sec:convergence}} is dedicated to the proof of an analytic
result on the convergence of nearly harmonic maps between surfaces to
a genuinely harmonic map with values in an $\R$-tree.  We take as
inspiration the work of Wolf on harmonic maps to surfaces and
$\R$-trees, extending his results for eventual application to the
collapsing map of a grafted surface.

\sboldpoint{Section \ref{sec:mainproof}} completes the proof of
Theorem \ref{thm:limit} by considering the collapsing map of a
divergent sequence of prunings, and applying the convergence result
from §\ref{sec:convergence}.  The same type of argument also allows us
to bound the difference between the Hopf and Hubbard-Masur
differentials associated to a grafted surface (Theorem
\ref{thm:hopfstrebel}), which can be seen as a finite analogue of the
asymptotic statement in Theorem \ref{thm:limit}.

\sboldpoint{Appendix \ref{sec:teichmuller}} contains a brief discussion of
Teichmüller geodesics in relation to the main results of the paper.
In particular we exhibit an asymmetry phenomenon that precludes the
existence of an extension of the Teichmüller geodesic involution to
the Thurston boundary of Teichmüller space.

\boldpoint{Acknowledgements.}  The results in this paper are taken
from my PhD thesis \cite{Du}, completed at Harvard University.  I
am grateful to my advisor, Curt McMullen, for his advice and
inspiration.  I would also like to thank the referee for a number of
insightful comments and helpful suggestions.

\section{Grafting, pruning, and collapsing}

\label{sec:grafting}

Let $S$ be a compact oriented surface of genus $g>1$, and $\T(S)$ the
Teichmüller space of marked conformal (equivalently, hyperbolic)
structures on $S$.  The simple closed hyperbolic geodesics on any
hyperbolic surface $Y \in \T(S)$ are in one-to-one correspondence with
the free homotopy classes of simple closed curves on $S$; therefore,
when a particular hyperbolic metric is under consideration, we will
use these objects interchangeably.

\begin{figure}
\begin{center}
\includegraphics[width=12cm]{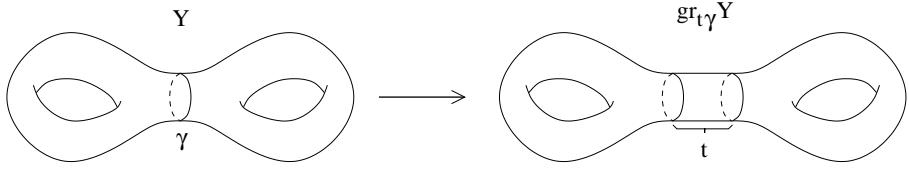}
\caption{The basic example of grafting.}
\label{fig:basicgraft}
\end{center}
\end{figure}

Fix $Y \in \T(S)$ and $\gamma$, a simple closed hyperbolic geodesic on $Y$.
\emph{Grafting} is the operation of removing $\gamma$ from $Y$ and
replacing it with a Euclidean cylinder $\gamma × [0,t]$, as shown in Figure
\ref{fig:basicgraft}.  The resulting surface is called the grafting of
$Y$ along the weighted geodesic $t \gamma$, written $\gr_{t \gamma} Y$.

Associated to each grafted surface $\gr_{t \gamma} Y$ is a canonical map
$\cl : \gr_{t \gamma} Y \to Y$, the \emph{collapsing map}, that collapses the
grafted cylinder $\gamma× [0,t]$ back onto the geodesic $\gamma$.  There is also
a natural $C^1$ conformal metric $\rhoth$ on $\gr_{t \gamma}Y$, the
\emph{Thurston metric}, that unites the hyperbolic metric on $Y$ with
the Euclidean metric of the cylinder $\gamma × [0,t]$.  The collapsing map
is distance nonincreasing with respect to the hyperbolic metric on
$Y$ and the Thurston metric on $\gr_{t \gamma}Y$.

Grafting is compatible with the natural generalization of weighted
simple closed geodesics to measured geodesic laminations.  The space
$\ML(S)$ of measured geodesic laminations is a contractible
PL-manifold with an action of $\R^+$ in which the set of weighted simple
closed hyperbolic geodesics for any hyperbolic metric $Y \in \T(S)$
forms a dense set of rays.  Conversely, given $Y \in \T(S)$ each
lamination $\lambda \in \ML(S)$ corresponds to a foliation of a closed subset of $Y$
by complete nonintersecting hyperbolic geodesics equipped with a
transverse measure of full support.  For a detailed treatment of
measured laminations, we refer the reader to \cite{Th2}.

Thurston has shown that grafting along simple closed curves extends
continuously to arbitrary measured laminations, and thus defines a map
$$ \gr : \ML(S) × \T(S) \to \T(S) \; \text{ where } \; (\lambda,Y) \mapsto \gr_\lambda
Y.$$
Morally, $\gr_\lambda Y$ is obtained from $Y$ by thickening the geodesic
realization of $\lambda$ in a manner determined by the transverse measure.

As in the simple closed curve case, there is a collapsing map $\cl :
\gr_\lambda Y \to Y$ that collapses the grafted part $A \subset \gr_\lambda Y$ onto the
geodesic realization of $\lambda$ on $Y$, and a conformal metric $\rhoth$ on
$\gr_\lambda Y$ that is hyperbolic on $\gr_\lambda Y - A$.  These can be defined
simply as limits of the corresponding objects for grafting along
simple closed curves, or in terms of the geometry of the locally
convex pleated plane in $\H^3$ associated to $(Y,\lambda)$ \cite{KT}.  The
metric $\rhoth$ is of class $C^{1,1}$ on $\gr_\lambda Y$, and thus its
curvature is defined almost everywhere \cite{KP2}.

Scannell and Wolf have shown that for each $\lambda \in \ML(S)$, the map
$\gr_\lambda: \T(S) \to \T(S)$ is a homeomorphism \cite{SW}.  Thus there is the inverse,
or \emph{pruning} map
$$ \pr : \ML(S) × \T(S) \to \T(S) \; \text{ where } \; (\lambda,X) \mapsto \pr_\lambda X =
\gr_\lambda^{-1}X.$$ 
In other words, for each $X \in \T(S)$ and $\lambda \in \ML(S)$,
there is a unique way to present $X$ as a grafting of some Riemann
surface $Y = \gr_\lambda^{-1}(X)$ along $\lambda$, and pruning is the operation of
recovering $Y$ from the pair $(\lambda,X)$.

We will be primarily interested in the pruning map when the surface
$X$ is fixed, i.e. the map $\lambda \mapsto \pr_\lambda X$ from
$\ML(S)$ to $\T(S)$.  Theorem \ref{thm:limit} describes the asymptotic
behavior of this map in terms of the conformal geometry of $X$.  Since
the graph of this map is exactly the preimage of $X$ under the
grafting map, the continuous extension of the former to a
homeomorphism of $\PML(S)$ also describes the closure of a fiber of
$\gr : \ML(S) × \T(S) \to \T(S)$; this is the content of
Theorem \ref{thm:fiber}.

Grafting and pruning are also related to $\CP^1$-structures on
surfaces, i.e. geometric structures modeled on $(\CP^1, \PSL_2(\C))$.
Each $\CP^1$ structure has an underlying conformal structure, and thus
the deformation space $¶(S)$ of marked $\CP^1$-surfaces on $S$ has a
natural projection map $\pi : ¶(S) \to \T(S)$ to Teichmüller space.

Thurston introduced an extension of grafting that produces complex
projective structures (see \cite{KT}, \cite{Ta}), and showed that the resulting map
$$ \Gr : \ML(S) × \T(S) \to ¶(S)$$
is a homeomorphism.  The composition of this projective grafting map
with the projection $\pi$ yields the conformal grafting map discussed
above, i.e. $\pi \circ \Gr = \gr$.  Thus the fiber $P(X) = \pi^{-1}(X)
\subset ¶(S)$ over $X \in \T(S)$ corresponds via $\Gr$ to the fiber
of the conformal grafting map:
$$ \Gr^{-1}(P(X)) = \gr^{-1}(X). $$
Applying this observation we see that Theorem \ref{thm:boundary} is a
restatement of Theorem \ref{thm:fiber} in the language of $\CP^1$
geometry.

\section{Conformal metrics and quadratic differentials}

\label{sec:tensors}

Our study of the pruning map will center on the collapsing map as a
candidate for the harmonic map variational problem (i.e. minimizing
energy in a given homotopy class). In this section we collect
definitions and background material required for this variational
analysis.

Fix a hyperbolic Riemann surface $X$ and
let $S(X)$ denote the space of measurable complex-valued quadratic
forms on $TX$.  A form $\beta \in S(X)$
can be decomposed according to the complex structure on $X$, i.e. 
$$ \beta = \beta^{2,0} + \beta^{1,1} + \beta^{0,2} $$ where $\beta^{i,j}$ is a symmetric
form of type $(i,j)$.  Then $|\beta^{1,1}|$ is a conformal metric on $X$,
which can be thought of as a ``circular average'' of $|\beta|$,
$$ |\beta^{1,1}(v)| = \frac{1}{2 \pi} \int_0^{2 \pi} |\beta(R_\theta v, R_\theta v)| d\theta,$$
where $v \in T_x X$ and $R_\theta \in \aut(T_x X)$ is the rotation by angle
$\theta$ defined by the conformal structure of $X$.
The area of $X$ with respect to $|\beta^{1,1}|$, when it is finite, defines
a natural $L^1$ norm
$$ \| \beta \|_{L^1(X)} = \int_X |\beta^{1,1}|, $$ 
which we will abbreviate to $\|\beta\|_1$ if the domain $X$ is fixed.  We
also call $\|\beta\|_1$ the \emph{energy} of $\beta$.

We write $S^{2,0}(X)$ for the space of \emph{measurable quadratic
  differentials} on $X$, which are those quadratic forms $\phi \in S(X)$
  such that $\phi = \phi^{2,0}$.  Within $S^{2,0}(X)$, there is the
  space $Q(X)$ of \emph{holomorphic quadratic differentials},
  i.e. holomorphic sections of $K_X^2$.  By the Riemann-Roch
  theorem, if $X$ is compact and has genus $g$, $ \dim_\C Q(X) = 3g-3$.

Fixing a conformal metric $\sigma$ also allows us to define the unit tangent
bundle
$$UTX = \left \{ (x,v) \in TX \: | \; \| v \|_\sigma = 1
\right \},$$
in which case the interpretation of $\beta^{1,1}$ as the circular average
of $\beta$ yields another expression for $\|\beta\|_{L^1(X)}$:
$$ \| \beta \|_{L^1(X)} = \frac{1}{2 \pi} \int_X \int_{UT_x X} |\beta(v)| d\theta(v) d\sigma(x) $$

The \emph{Hopf differential} $\Phi(\beta)$ of $\beta \in S(X)$ is the $(2,0)$
part of its decomposition,
$$ \Phi(\beta) = \beta^{2,0},$$
which (along with $\Phi(\overline{\beta})$) measures the failure of
$\beta$ to be compatible with the conformal structure of $X$.  For
example, there is a function $b : X \to \C$ such that $\beta = b \sigma$ if
and only if $\Phi(\beta) = \Phi(\overline{\beta}) = 0$.

Let $f : X \to (M,\rho)$ be a smooth map from $X$ to a Riemannian manifold
$(M,\rho)$.  Then the energy $\E(f)$ and Hopf differential $\Phi(f)$ of $f$
are defined to be those of the pullback metric $f^*(\rho)$:
\begin{equation}
\label{eqn:energyhopf}
\begin{split}
\E(f) &= \| f^*(\rho) \|_1 =  \int_X |f^*(\rho)^{1,1}|\\
\Phi(f) &= [f^*(\rho)]^{2,0}
\end{split}
\end{equation}
Since $f^*(\rho)$ is a real quadratic form, at points where $Df$ is
nondegenerate $f$ is conformal if and only if $\Phi(f) = 0$.  Thus
$\E(f)$ is a measure of the average stretching of the map $f$, while
$\Phi(f)$ records its anisotropy.

\section{Measured foliations and the antipodal map}

\label{sec:foliations}

In this section we briefly recall the identifications between the
spaces of measured foliations, measured geodesic laminations, and
holomorphic quadratic differentials on a compact Riemann surface.  We
then use these identifications to define the antipodal map on the
space of measured laminations.  

A measured foliation on $S$ is a singular foliation $\F$ (i.e. one
with isolated $k$-pronged singularities) and an assignment of a Borel
measure to each transverse arc in a manner compatible with
transversality-preserving isotopy.  A detailed discussion of measured
foliations can be found in \cite{FLP}.  The notation $\MF(S)$ is used
for the quotient of the set of measured foliations by the equivalence
relation generated by isotopy and Whitehead moves (e.g., collapsing
leaves connecting singularities).

The typical example of a measured foliation comes from a holomorphic
quadratic differential $\phi \in Q(X)$.  The measured foliation
$\F(\phi)$ determined by $\phi$ is the pullback of the horizontal line
foliation of $\C$ under integration of the locally defined holomorphic
$1$-form $\sqrt{\phi}$.  Equivalently, a vector $v \in T_xX$ is
tangent to $\F(\phi)$ if and only if $\phi(v) > 0$.  The measure on
transversals is obtained by integrating the length element $|\im
\sqrt{\phi}|$.

The foliation $\F(\phi)$ is called the \emph{horizontal foliation} of
$\phi$.  Since $\sqrt{-\phi} = i \sqrt{\phi}$, $\F(\phi)$ and $\F(-\phi)$ are
orthogonal, and $\F(-\phi)$ is called the \emph{vertical foliation} of
$\phi$.

Hubbard and Masur proved that measured foliations and quadratic
differentials are essentially equivalent notions:

\begin{thm}[Hubbard and Masur, \cite{HM}]
For each $\nu \in \MF(S)$ of measured foliations and $X \in \T(S)$ there is
a unique holomorphic quadratic differential $\phi_X(\nu) \in Q(X)$ such that
$$ \F(\phi_X(\nu)) = \nu.$$ Furthermore, the map $\phi_X : \MF(S) \to Q(X)$ is a
homeomorphism.
\end{thm}

Note that the transverse measure of $\phi \in Q(X)$ is defined using
$\sqrt{\phi}$, and so for $c > 0$,
$$ \F(c \phi) = c^{\frac{1}{2}} \F(\phi).$$ As a result, the Hubbard-Masur
map $\phi_X$ has the following homogeneity property:
\begin{equation}
\label{eqn:hmquadratic}
\phi_X(C \nu) = C^2 \phi_X(\nu) \text{ for all } C >0
\end{equation} 

One can also view measured foliations as diffuse versions of measured
laminations, in that every measured foliation is associated to a
unique measured lamination with the same intersection properties (for
details see \cite{Lev}).  This induces a homeomorphism between
$\MF(S)$ and $\ML(S)$.  Using this homeomorphism implicitly, we can
consider the Hubbard-Masur map $\phi_X$ to have domain $\ML(S)$, and
we write $\Lambda$ for its inverse,
$$ \Lambda : Q(X) \to \ML(S).$$

We can also use $\Lambda$ to transport the linear involution $\phi \mapsto (-\phi)$ of
$Q(X)$ to an involutive homeomorphism $i_X : \ML(S) \to \ML(S)$, i.e.
$$ i_X(\lambda) = \Lambda(-\phi_X(\lambda)).$$ 
Since $\F(\phi)$ and $\F(-\phi)$ are orthogonal
foliations, we say that $\lambda, \mu \in \ML(S)$ are \emph{orthogonal with
respect to $X$} if $i_X(\lambda) = \mu$.

The resulting homeomorphism depends on $X \in \T(S)$ in an essential
way, just as the orthogonality of foliations or tangent vectors
depends on the choice of a conformal structure.

Since $\Lambda$ is homogeneous, it also induces a homeomorphism between
projective spaces:
$$ \Lambda : \PQ(X) = (Q(X)-\{0\}) / \R^+ \to \PML(S) = (\ML(S)-\{0\}) / \R^+.$$
Thus we obtain an involution $i_X : \PML(S) \to
\PML(S)$ that is topologically conjugate to the antipodal map $(-1) :
\PQ(X) \simeq S^{2n - 1} \to S^{2n-1}$.  We call $i_X$ the \emph{antipodal
  involution with respect to $X$}.

\section{$\R$-trees from measured laminations and foliations}

\label{sec:rtrees}

An $\R$-tree (or \emph{real tree}) is a complete geodesic metric space
in which there is a unique embedded path joining every pair of points,
and each such path is isometric to an interval in $\R$.  Such
$\R$-trees arise naturally in the context of measured foliations on
surfaces, as we now describe (for a detailed treatment, see \cite{Ka}).

Let $\F \in \MF(S)$ be a measured foliation on $S$, and lift $\F$ to a
measured foliation $\Tilde{\F}$ of the universal cover $\Tilde{S}
\simeq \H^2$.  Define a pseudometric $d_\F$ on $\Tilde{S}$,
\begin{equation}
\label{eqn:pseudometric}
 d_\F(x,y) = \inf \: \{ i(\Tilde{\F},\gamma) \: | \: \gamma : [0,1] \to \Tilde{S},\:
  \gamma(0)=x,\: \gamma(1) = y \}
\end{equation}
where $i(\Tilde{\F},\gamma)$ is the intersection number of a transverse arc
$\gamma$ with the measured foliation $\Tilde{\F}$.

The quotient metric space $T_\F = \Tilde{S} / ( x \sim y \:
  \text{ if } \: d_\F(x,y) = 0)$ is an $\R$-tree whose isometry type
  depends only on the measure equivalence class of $\F$.  Alternately,
  $T_\F$ is the space of leaves of $\Tilde{\F}$ with metric induced by
  the transverse measure, where we consider all of the leaves that
  emanate from a singular point of $\Tilde{\F}$ to be a single point
  of $T_\F$.  The action of $\pi_1(S)$ on $\Tilde{S}$ by deck
  transformations descends to an action on $T_\F$ by isometries; in
  fact, by Skora's theorem, every such isometric action of $\pi_1(S)$ on
  an $\R$-tree satisfying a nondegeneracy condition (\emph{small} and
  \emph{minimal}) arises in this way (see \cite{Skora}, also \cite{FW}, \cite{Otal}).

If $\F$ is a measured foliation associated to the measured lamination
  $\lambda$, then the resulting $\R$-tree $T_\lambda$ with metric $d_\lambda$ is called
  the \emph{dual $\R$-tree} of $\lambda$.  For example, when $\lambda$ is
  supported on a family of simple closed curves, the $\R$-tree
  $T_\lambda$ is actually a simplicial tree of infinite valence with one
  vertex for each lift of a complementary region of $\lambda$ to
  $\Tilde{S}$, and where an edge connecting two adjacent complementary
  regions has length equal to the weight of the geodesic that
  separates them.

A slight generalization of this construction arises naturally in the
context of grafting.  The grafting locus $A = \cl^{-1}(\lambda) \subset \gr_\lambda Y$ has a natural
foliation $\F_A$ by Euclidean geodesics that map isometrically onto
$\lambda$, with a transverse measure induced by the Euclidean metric in the
orthogonal direction. The associated pseudometric on $\widetilde{\gr_\lambda
Y}$,
$$ d_{\F_A}(x,y) = \inf \: \{ i(\Tilde{\F_A},\gamma) \: | \: \gamma : [0,1] \to
  \widetilde{\gr_\lambda Y},\:
  \gamma(0)=x,\: \gamma(1) = y \},$$
yields a quotient $\R$-tree isometric to $T_\lambda$ and a map
$$ \hat{\cl} : \widetilde{\gr_\lambda Y} \to T_\lambda, $$ which we
call the \emph{co-collapsing map}.  While the collapsing map $\cl :
\gr_\lambda Y \to Y$ compresses the entire grafted part back to its
geodesic representative, the co-collapsing map collapses each
connected component of $(\widetilde{\gr_\lambda Y} - \tilde{A})$, i.e
the complement of the grafted part of $\widetilde{\gr_\lambda Y}$, and
each leaf of $\F_A$ to a single point.

\section{Harmonic maps to surfaces and $\R$-trees}

\label{sec:harm}

We now consider the harmonic maps variational problem, i.e. finding
critical points of the energy functional $\E(f)$ for maps $f : X \to
(M,\rho)$ from a Riemann surface $X$ to a Riemannian manifold $M$.  In
§\ref{sec:tensors}, we defined the energy for smooth maps, but the
natural setting in which to work with the energy functional is the
Sobolev space $W^{1,2}(X,M)$ of maps with $L^2$ distributional
derivatives (see \cite{KS} for details).  The maps we consider are
Lipschitz, hence differentiable almost everywhere, so there is no
ambiguity in our definitions.

A critical point of the energy functional $\E : W^{1,2}(X,M) \to \R^{\geq
  0}$ is a \emph{harmonic map}; it is easy to show that the Hopf differential $$
  \Phi(f) = f^*(\rho)^{2,0} $$ is holomorphic ($\Phi(f) \in Q(X)$) if $f$ is
  harmonic.  In particular, the maximal and minimal stretch directions
  of a harmonic map are realized as a pair of orthogonal foliations by
  straight lines in the singular Euclidean metric $|\Phi(f)|$.

For any pair of compact hyperbolic surfaces $X,Y$ and nontrivial
homotopy class of maps $[f]: X \to Y$, there is a unique harmonic map $h
: X \to Y$ that is smooth and homotopic to $f$; furthermore, $h$
minimizes energy among such maps \cite{Ha}.  In particular, for each
$X,Y \in \T(S)$ there is a unique harmonic map $h : X \to Y$ that is
compatible with the markings; define
$$ \E(X,Y) = \E(h : X \to Y).$$

For the proof of Theorem \ref{thm:limit} we will also need to consider
harmonic maps from Riemann surfaces to $\R$-trees.  The main
references for this theory are the papers of Wolf (\cite{W4},
\cite{W5}, \cite{W6}); a much more general theory of harmonic maps to
metric spaces is discussed by Korevaar and Schoen in \cite{KS} and
\cite{KS2}.  The $\R$-trees and harmonic maps we consider will arise
from limits of maps to degenerating hyperbolic surfaces; with the
appropriate notion of convergence, Wolf has shown that an $\R$-tree
limit can always be extracted:

\begin{thm}[Wolf \cite{W4}]
\label{thm:rtreelimit}
Let $Y_i \in \T(S)$ and $\mu \in \ML(S)$ be such that $Y_i \to [\mu] \in \PML(S)$
in the Thurston compactification.  Then after rescaling the hyperbolic
metrics $\rho_i$ on $Y_i$ appropriately, the sequence of universal
covering metric spaces $(\Tilde{Y_i}, \Tilde{\rho_i})$
converges in the equivariant Gromov-Hausdorff sense to the $\R$-tree
$T_\mu$.
\end{thm}

The equivariant Gromov-Hausdorff topology is a natural setting in
which to consider convergence of metric spaces equipped with isometric
group actions.  The application of this topology to the Thurston
compactification is due to Paulin \cite{Pa}; for other perspectives
on the connection between $\R$-trees and Teichmüller theory, see
\cite{Be}, \cite{MS}, \cite{CS}.

The theory of harmonic maps is well-adapted to this
generalization from smooth surfaces to metric spaces like $\R$-trees;
for example, the convergence statement of Theorem \ref{thm:rtreelimit}
can be extended to a family of harmonic maps from a fixed Riemann
surface $X$:

\begin{thm}[Wolf \cite{W4}]
\label{thm:harmconverge}
For $X \in \T(S)$ and $\lambda \in \ML(S)$, let $\pi_\lambda : \Tilde{X} \to T_\lambda$
denote the projection onto the leaves of $\F(\phi_X(\lambda))$, where $\phi_X(\lambda)$
is the Hubbard-Masur differential for $\lambda$.  Then:
\begin{enumerate}
\item  $\pi_\lambda$ is harmonic, meaning that it pulls back germs of convex functions on
  $T_\lambda$ to germs of subharmonic functions on $\Tilde{X}$.
\item If $Y_i \in \T(S)$ is a sequence such that $Y_i \to [\lambda]
  \in \PML(S)$, then the lifts $\Tilde{h_i} : \Tilde{X} \to
  \Tilde{Y_i}$ of the harmonic maps converge in the equivariant Gromov-Hausdorff sense
  to $\pi_X : \Tilde{X} \to T_{\lambda}$.
\end{enumerate}
\end{thm}

Much like the case of maps between Riemann surfaces, it is most
natural to work with the energy functional on a Sobolev space
$W^{1,2}(\Tilde{X}, T_\lambda)$ of equivariant maps with $L^2$
distributional derivatives, and again we defer to \cite{KS} for details,
since the maps we consider are Lipschitz.  

By definition, the metric $d_\lambda$ on $T_\lambda$ is isometric to
the standard Euclidean metric $dx^2$ on $\R$ along each geodesic
segment.  Thus if $f : \Tilde{X} \to T_\lambda$ is an equivariant
Lipschitz map that sends a neighborhood of almost every point in
$\Tilde{X}$ into a geodesic segment of $T_\lambda$, then the pullback metric
$f^*(dx^2)$ is a well-defined (possibly degenerate) measurable
quadratic form on $T\Tilde{X}$ which is invariant under the action of
$\pi_1(X)$.  This allows us to define the energy $\E(f)$ and Hopf
differential $\Phi(f)$ of such an equivariant map as the $L^1$ norm
and $(2,0)$ part of the induced quadratic form on $TX$.

For later use we record the following calculations relating the Hopf
differential and energy of the projection $\pi_\lambda$ and the
Hubbard-Masur differential $\phi_X(\lambda)$; details can be found in
\cite{W4}.

\begin{lem}
\label{lem:hopfenergyprojection}
The Hopf differential of $\pi_\lambda : \Tilde{X} \to T_\lambda$ is
$$ \Phi(\pi_\lambda) = - \phi_X(\frac{1}{2}\lambda) = - \frac{1}{4}\phi_X(\lambda),$$
and the energy of $\pi_\lambda$ is given by
$$ \E(\pi_\lambda) = \frac{1}{2} E(\lambda,X) = \frac{1}{2} \|
\phi_X(\lambda) \|_1, $$
where $E(\lambda,X)$ is the extremal length of $\lambda$ on $X$.
\end{lem}

\section{Energy and grafting}

\label{sec:energy}

In this section we recall Tanigawa's computation of the energy of the
collapsing map and discuss its consequences for grafting and pruning.
Tanigawa used this energy computation to show that for each $\lambda
\in \ML(S)$, grafting $\gr_\lambda : \T(S) \to \T(S)$ is a
proper map of Teichmüller space.

\begin{lem}[Tanigawa \cite{Ta}]
\label{lem:tanigawa2}
Let $X = \gr_\lambda Y$, where $X,Y \in \T(S)$ and $\lambda \in \ML(S)$, and let $h$
denote the harmonic map $h : X \to Y$ compatible with the markings
(and thus homotopic to $\cl : X \to Y$).  Then

$$ \frac{1}{2} \ell(\lambda,Y) \leq \frac{1}{2} \frac{\ell(\lambda,Y)^2}{E(\lambda,X)} \leq \E(h) \leq
\frac{1}{2} \ell(\lambda,Y) + 2 \pi |\chi(S)| = \E(\cl),$$
where $\ell(\lambda,Y)$ is the hyperbolic length of $\lambda$ on $Y$, and
$E(\lambda,X)$ is the extremal length of $\lambda$ on $X$.
\end{lem}

The middle part of Tanigawa's inequality, i.e.
$$ \frac{1}{2} \frac{\ell(\lambda,Y)^2}{E(\lambda,X)} \leq \E(h) $$ is due to Minsky,
and holds for any harmonic map between finite-area Riemann surfaces of
the same type and any measured lamination $\lambda$ \cite{Mi}.

We will also need to know the energy of the co-collapsing map, which
is easily computed using the same method as in \cite{Ta}:

\begin{lem}
\label{lem:cocollapseenergy}
Let $X = \gr_\lambda Y$, where $X,Y \in \T(S)$ and $\lambda \in
\ML(S)$.  Then the energy $\E(\Hat{\cl})$ of the co-collapsing map
$\Hat{\cl} : \Tilde{X} \to T_\lambda$ is given by
$$ \E(\Hat{\cl}) = \frac{1}{2} \ell(\lambda,Y).$$
\end{lem}

\begin{proof}
On the grafting locus, both the collapsing and co-collapsing maps are
locally modeled on the orthogonal projection of a Euclidean plane onto
a geodesic.  Thus the conformal parts of the pullback metrics
associated to the collapsing and co-collapsing maps are identical on
the grafting locus, i.e.
$$ [\cl^*(\rho)]^{1,1} = [\Hat{\cl}^*(\rho)]^{1,1}.$$ On the remainder
of the surface (the hyperbolic part), the pullback metric via the
co-collapsing map is zero, while the collapsing map is an isometry
here.  The energy difference is therefore the hyperbolic area of $Y$,
$$ \E(\cl) = \E(\Hat{\cl}) + 2 \pi |\chi(S)|,$$
and the lemma follows.
\end{proof}

It is an immediate consequence of Lemma \ref{lem:tanigawa2} that the
the collapsing map has energy very close to that of the harmonic map:

\begin{cor}
\label{cor:energyclose}
Let $h:\gr_\lambda Y \to Y$ be the harmonic map homotopic to
the collapsing map $\cl$.  Then
$$\E(h) < \E(\cl) \leq \E(h) + 2 \pi |\chi(S)|.$$
Note that the constant on the right hand side is independent of $Y \in
\T(S)$ and $\lambda \in \ML(S)$.
\end{cor}

For our purposes, another important consequences of Lemma
\ref{lem:tanigawa2} is a relationship between the hyperbolic length of the
grafting lamination on $Y$ and its extremal length on $\gr_\lambda Y$:

\begin{lem}
\label{lem:lengthextremal}
Let $X = \gr_\lambda Y$, where $Y \in \T(S)$.  Then we have
$$ \ell(\lambda, Y) =  E(\lambda,X) + O(1) $$
where the implicit constant depends only on $|\chi(S)|$.
\end{lem}

\begin{proof}
The lower bound
$$ \ell(\lambda, Y) \geq E(\lambda, X)$$
is immediate from the left-hand side of the inequality of Lemma
\ref{lem:tanigawa2}, and we also have
$$\frac{1}{2} \frac{\ell(\lambda,Y)^2}{E(\lambda,X)}  \leq \frac{1}{2} \ell(\lambda,Y) + 2 \pi
|\chi(S)|$$
and therefore,
$$ \frac{\ell(\lambda,Y)^2}{\frac{1}{2} \ell(\lambda,Y) + O(1)} \leq E(\lambda,X). $$
Solving for $\ell(\lambda, Y)$ yields
\begin{equation*}
\begin{split}
\ell(\lambda, Y) &\leq \frac{1}{2} E(\lambda,X) + \frac{1}{2} (E(\lambda,X))^{1/2} (E(\lambda,
X) + O(1))^{1/2}\\
& = E(\lambda, X) + O(1)
\end{split}
\end{equation*}
\end{proof}

Using Lemma \ref{lem:cocollapseenergy} and this comparison between
hyperbolic and extremal length, we find that the energy of the
co-collapsing map is close to that of the harmonic projection
(cf. Corollary \ref{cor:energyclose}):
\begin{cor}
\label{cor:coenergyclose}
Let $\pi_\lambda$ be the harmonic projection from $\Tilde{X} =
\widetilde{\gr_\lambda Y}$ to the $\R$-tree $T_\lambda$, and
$\Hat{\cl}$ the associated co-collapsing map.  Then
$$\E(\pi_\lambda) < \E(\Hat{\cl}) \leq \E(\pi_\lambda) + 2 \pi |\chi(S)|.$$
\end{cor}

Finally, as a complement to Tanigawa's result on grafting, we can use
these energy computations to show that the pruning map with basepoint
$X$ is proper:

\begin{lem}
\label{lem:divergence}
For each $X \in \T(S)$, the pruning map with basepoint $X$, i.e.
$$ \pr_{(\cdot)} X : \ML(S) \to \T(S)$$
is proper.
\end{lem}

\begin{proof}
Suppose on the contrary that $\lambda_i \to \infty$ but the
sequence $Y_i = \pr_{\lambda_i} X$ remains in a compact subset of
Teichmüller space. Then $\ell(\lambda_i, Y_i) \to \infty$ and by
Lemma \ref{lem:tanigawa2},
$$ \E(h_i) \geq \frac{1}{2} \ell(\lambda_i,Y_i) \to \infty,$$
where $h_i : X \to Y_i$ is the harmonic map compatible with
the markings.

On the other hand, a result of Wolf (see \cite{W}) states that for any
fixed $X \in \T(S)$, the energy $\E(X, \cdot )$ is a proper function on
$\T(S)$.  Since $\E(h_i) \to \infty$, we conclude $Y_i \to \infty$, which is a
contradiction.
\end{proof}

\section{Hopf differentials and grafting}

\label{sec:hopf}

We will need to consider not only the energy but also the Hopf
differentials $\Phi(\cl)$ and $\Phi(\Hat{\cl})$ associated to the
collapsing and co-collapsing maps of a grafted surface.  In
particular, for the proof of Theorem \ref{thm:limit}, we will use a
relationship between these differentials and the grafting lamination.
In this section we establish such a relationship, after addressing some
technical issues that arise because the quadratic differentials under
consideration are not holomorphic.

For holomorphic quadratic differentials $\phi, \psi \in Q(X)$, the
intersection number of their measured foliations can be expressed in
terms of the differentials (see \cite{Ga}); define
$$ \omega(\phi, \psi)= \int_X | \im \left ( \sqrt{\phi} \sqrt{\psi} \right ) |.$$
Then
$$ i(\F(\phi), \F(\psi)) = \omega(\phi, \psi).$$ However, the quantity $\omega(\alpha, \beta)$ makes
sense for $L^1$ quadratic differentials $\alpha$ and $\beta$, holomorphic or
not.

While a measurable differential $\alpha$ does not define a measured
foliation, it does have a \emph{horizontal line field}
$\mathscr{L}(\alpha)$ consisting of directions $v \in TX$ such that $\alpha(v) >
0$.  Then $\omega(\alpha,\beta)$ measures the average transversality (sine of twice
the angle) between the line fields $\mathscr{L}(\alpha)$ and
$\mathscr{L}(\beta)$, averaged with respect to the measure $|\alpha|^{1/2}
|\beta|^{1/2}$.

Now consider the collapsing map $\cl : X \to \pr_\lambda X$, and for
simplicity let us first suppose $\lambda$ is supported on a single simple
closed hyperbolic geodesic $\gamma \subset \pr_\lambda X$, i.e. $\lambda = t \gamma$.  Then the
grafting locus $A \subset X$ is the Euclidean cylinder $\gamma × [0,t]$, and
the collapsing map is the projection onto the geodesic $\gamma$.  Just
as the Hopf differential of the orthogonal projection of $\C$ onto
$\R$ is
$$ \Phi(z \mapsto \re(z)) = [ dx^2 ]^{2,0} = \frac{1}{4} dz^2,$$
the Hopf differential of $\cl$ on $A$ is the pullback of
$\frac{1}{4}dz^2$ via local Euclidean charts that take parallels of
$\gamma$ to horizontal lines.  This differential is holomorphic on
$A$, and corresponds to the measured foliation $\frac{1}{2} \lambda$.
Thus the Euclidean metric on $A$, which is the restriction of the Thurston
metric of $X$, is given by $| 4 \Phi(\cl) |$.

On the complement of the grafting locus, the collapsing map is conformal
and thus the Hopf differential is zero.  Therefore $\Phi(\cl)$ is a
piecewise holomorphic differential on $X$ whose horizontal line
field is the natural foliation of the grafting locus by parallels of
the grafting lamination, with half of the measure of $\lambda$.  This
analysis extends by continuity to the case of a general lamination $\lambda
\in \ML(S)$.

It follows that the line field $\mathscr{L}(\Phi(\cl))$ represents the
measured lamination $\frac{1}{2}\lambda$, in that for all $\psi \in Q(X)$,
\begin{equation}
\label{eqn:intersection} 
\omega(\Phi(\cl), \psi) = \frac{1}{2} i(\lambda,\F(\psi)).
\end{equation}
We therefore use the notation
$$ \Phi_X(\lambda) = \Phi( \cl : X \to \pr_\lambda X)$$
for the Hopf differential of the collapsing map, which is somewhat like $\phi_X(\frac{1}{2}\lambda)$
in that it is a quadratic differential whose foliation is a 
distinguished representative for the measured foliation class of
$\frac{1}{2}\lambda$.  The Hopf differential $\Phi_X(\lambda)$ is not holomorphic, however,
though we will later see (§\ref{sec:mainproof}) that is is nearly so.

For now, we will simply show that $L^1$ convergence of Hopf
differentials $\Phi_X(\lambda)$ to a holomorphic limit implies convergence
of the laminations $\lambda$:

\begin{lem}
\label{lem:intersection}
Let $X \in \T(S)$ and $\lambda_i \in \ML(S)$.  If
$$ [\Phi_X(\lambda_i)] \to [\psi], \; \text{ where } \psi \in Q(X)$$
then
$$ [\lambda_i] \to [\F(\psi)] \in \PML(S).$$
Here $[\Phi_X(\lambda_i)]$ is the image of $\Phi_X(\lambda_i)$ in $\mathbb{P}S^{2,0}(X)$.
\end{lem}

\begin{proof}
First, we can choose $c_i > 0$ such that
$$ c_i^2 \Phi_X(\lambda_i) \to \psi.$$

It is well known that there are finitely many simple closed curves
$\nu_k$, $k = 1 \ldots N$, considered as measured laminations with
unit weight, such that the map $I : \ML(S) \to \R^N$ defined by
$I(\lambda) = \{i(\lambda,\nu_k)\}$ is a homeomorphism onto its image.
Recall that $\phi_X(\nu_k) \in Q(X)$ is the unique holomorphic
quadratic differential satisfying $\F(\phi_X(\nu_k)) = \nu_k$.

Since $\omega(\cdot,\nu_k) : S^{2,0}(X) \to \R$ is evidently a continuous map,
we conclude from (\ref{eqn:intersection}) and the hypothesis
$c_i \Phi_X(\lambda_i) \to \psi$ that
$$ \omega(c_i^2 \Phi_X(\lambda_i), \phi_X(\nu_k)) = \frac{c_i}{2} i(\lambda_i, \nu_k) \to \frac{1}{2}
i(\F(\psi), \nu_k).$$ and
therefore $c_i \lambda_i \to \frac{1}{2}\F(\psi)$.
\end{proof}

The connection between the Hopf differential $\Phi_X(\lambda)$ and the Thurston
metric also allows us to compute its norm and relate it to the extremal
length of $\lambda$:

\begin{cor}
\label{cor:normhopf}
The $L^1$ norm of $\Phi_X(\lambda)$ is given by
$$ \| \Phi_X(\lambda) \|_1 = \frac{1}{4} \ell(\lambda, \pr_\lambda X)
 = \frac{1}{4} E(\lambda, X) + O(1).$$
\end{cor}

\begin{proof}
We have seen that $|4 \Phi_X(\lambda)|$ induces the Thurston metric on the
grafted part $A \subset X$ and is zero elsewhere.  The area of $A$ with
respect to the Thurston metric is $\ell(\lambda, \pr_\lambda X)$, and therefore,
$$ \| \Phi_X(\lambda) \|_1 = \frac{1}{4}\ell(\lambda, \pr_\lambda
 X).$$ On the other hand, it follows from Lemma \ref{lem:lengthextremal}
 that
$$\ell(\lambda, \pr_\lambda X) = E(\lambda,X) + O(1).$$
\end{proof}

We can apply the same analysis to the Hopf differential of the
co-collapsing map $\hat{\cl} : \widetilde{\Gr_\lambda Y} \to T_\lambda$.
Though the co-collapsing map is defined on the universal cover of the
grafted surface, its Hopf differential is invariant under the action
of $\pi_1(S)$ and therefore descends to a measurable quadratic
differential $\Phi(\hat{\cl}) \in S^{2,0}(X)$.

The co-collapsing map $\hat{\cl}$ is piecewise constant on the
complement of the grafting locus, so (like $\Phi(\cl)$) its Hopf
differential is identically zero there.  Within the grafting locus it
is modeled on the orthogonal projection of $\C$ onto $i \R$ (where the
leaves of $\F_A$ correspond to horizontal lines in $\C$).  Since
$$ \Phi(z \mapsto \im(z)) = [ dy^2 ]^{2,0} = -\frac{1}{4} dz^2,$$ we conclude
that the Hopf differentials $\Phi(\hat{\cl})$ and $\Phi(\cl)$ are
inverses, i.e. 
\begin{equation}
\label{eqn:hopfhopf}
\Phi(\hat{\cl}) = - \Phi(\cl) = -\Phi_X(\lambda)
\end{equation}

\begin{rem}
The relationship between $\cl$ and $\hat{\cl}$ and their Hopf
differentials is reminiscent of the ``minimal suspension'' technique
introduced by Wolf; for details, see \cite{W7}.
\end{rem}

\section{Convergence to the harmonic map}

\label{sec:convergence}

We will use the fact that the collapsing map associated to a
grafted surface is nearly harmonic to extract a genuinely harmonic
limit map with values in an $\R$-tree.  In this section we formalize
the conditions on a sequence of maps that allows us to use this limit
construction and prove the main technical result about convergence.

Let $X, Y_i \in \T(S)$ and suppose $Y_i \to \infty$; let $\rho_i$ denote a
hyperbolic metric on $Y_i$, and $h_i : X \to Y_i$ the harmonic map (with
respect to $\rho_i$) compatible with the markings.

We say that a sequence of maps $f_i \in W^{1,2}(X,Y_i)$ compatible with
the markings of $X$ and $Y_i$ is a \emph{minimizing sequence} if
\begin{equation}
\label{eqn:minseq}
\lim_{i \to \infty} \frac{\E(f_i)}{\E(h_i)} = 1.
\end{equation}
Since the harmonic map $h_i$ is the unique energy minimizer in its
homotopy class, a minimizing sequence asymptotically minimizes energy.
In this section we will show that all minimizing sequences have the
same asymptotic behavior, in that their Hopf differentials have the
same projective limit.

Note that the condition that $f_i$ be a minimizing sequence implies
only that the energy differences $\E(f_i) - \E(h_i)$ are $o(\E(h_i))$,
while in the proof of Theorem \ref{thm:limit} we consider a sequence
of collapsing maps for which these energy differences are
$O(\E(h_i)^{\frac{1}{2}})$ (while $\E(h_i) \to \infty$), a stronger
condition.  We prove the next result for general minimizing sequences,
however, since a specific bound on the energy difference does not
substantially simplify its proof.

\begin{thm}
\label{thm:minimizing}
Let $X$ and $Y_i$ be as above, and suppose $\lim_{i \to \infty} Y_i = [\mu] \in
\PML(S)$ in the Thurston compactification.

Then for any minimizing sequence $f_i : X \to Y_i$, the measurable
quadratic differentials $[ f_i^*(\rho_i) ]^{2,0}$
converge projectively in the $L^1$ sense to a holomorphic quadratic
differential $\Phi \in Q(X)$ such that
$$ [ \F(-\Phi) ] = [ \mu ], $$ i.e. there are constants $c_i > 0$ such that
$$ \lim_{i \to \infty} c_i [f_i^*(\rho_i) ]^{2,0} = \Phi.$$
\end{thm}

\begin{note}
The vertical foliation $\F(-\Phi)$ appears in Theorem
\ref{thm:minimizing} because the Thurston limit is a lamination whose
intersection number provides an estimate of hyperbolic length;
directions orthogonal to $\F(-\Phi)$ (that is, tangent to $\F(\Phi)$) are
maximally stretched by a map with Hopf differential $\Phi$, so the
intersection number with $\F(-\Phi)$ provides such a length estimate.
\end{note}

Before giving the proof of Theorem \ref{thm:minimizing}, we recall a
theorem of Wolf upon which it is based.
\begin{thm}[Wolf, \cite{W}]
\label{thm:wolf}
Let $X, Y_i \in \T(S)$, and let $\Psi_i$ denote the Hopf differential of
the unique harmonic map $h_i:X \to Y_i$ respecting markings, where $Y_i$
is given the hyperbolic metric $\rho_i$.  Then
$$ \lim_{i\to\infty} Y_i = [\mu] \in \PML(S) $$ if and only if

$$ \lim_{i\to\infty} [\F(-\Psi_i)] = [\mu].$$
\end{thm}

In other words, if one compactifies Teichmüller space according to the
limiting behavior of the Hopf differential of the harmonic map from a
fixed conformal structure $X$, then the vertical foliation map $\F \circ (-1) : \PQ(X) \to
\PMF(S) \simeq \PML(S)$ identifies this compactification with the Thurston
compactification.

Theorem \ref{thm:wolf} is actually a consequence of the convergence of harmonic
maps $h_i$ to the harmonic projection $\pi_\mu$ to an $\R$-tree
(combining Theorem \ref{thm:rtreelimit} and Theorem
\ref{thm:harmconverge}), though in \cite{W} Wolf provides an
elementary and streamlined proof of this result.

We will compare the Hopf differentials of a minimizing sequence to
those of a harmonic map by means of the pullback metrics.  The
essential point is that the energy functional for maps to a
nonpositively curved (NPC) space is convex along geodesic homotopies.
Korevaar and Schoen use this to show that maps with nearly minimal
energy have nearly the same pullback metric:

\begin{blankthm}[{Korevaar-Schoen, \cite[\S2.6]{KS}}]
\label{thm:ksmain}
Let $X$ be a compact Riemann surface and $Z$ an NPC metric space
equipped with an isometric action of $\pi_1(X)$.  Let $f \in
W^{1,2}(\Tilde{X},Z)$ be an equivariant map and $h \in W^{1,2}(\Tilde{X},Z)$
an equivariant energy-minimizing map.  Then
\begin{equation}
\frac{1}{2 \pi} \int_{UTX} |\pi_f(v,v) - \pi_h(v,v)| \leq  \left ( 2
\E(f) - 2 \E(h) \right)^\frac{1}{2} \left ( \E(f)^\frac{1}{2} +
\E(h)^\frac{1}{2} \right )
\end{equation}
\end{blankthm}

Here $\pi_f$ is the quadratic form that plays the role of the
metric pulled back by $f$ in the Korevaar-Schoen theory of Sobolev
maps to NPC metric spaces (see \cite[\S2]{KS}).  Note that we have
specialized this theorem to the case of a compact surface domain,
whereas the actual estimate in \cite{KS} applies to any compact
Riemannian manifold.

Applying this theorem to the case where $Z$ the universal cover of a
compact hyperbolic surface or an $\R$-tree, in which case $\pi_f(v,v)
= f^*(\rho)(v)$ is the pullback metric considered above, we obtain the
following estimates:

\begin{cor}
\label{cor:energydiff}
Let $f \in W^{1,2}(X,Y)$ where $X,Y \in \T(S)$ and $Y$ is given the
hyperbolic metric $\rho$.  Let $h$ be the harmonic map homotopic to
$f$.  Then
\begin{equation*}
\| f^*(\rho) - h^*(\rho) \|_1 \leq \sqrt{2} \left ( \E(f) - \E(h)
\right )^\frac{1}{2} \left ( \E(f)^\frac{1}{2} + \E(h)^\frac{1}{2}
\right )
\end{equation*}
and in particular
\begin{equation*}
\| \Phi(f) - \Phi(h) \|_1 \leq \sqrt{2} \left ( \E(f) - \E(h)
\right )^\frac{1}{2} \left ( \E(f)^\frac{1}{2} + \E(h)^\frac{1}{2}
\right ).
\end{equation*}
\end{cor}

\begin{cor}
\label{cor:treeenergydiff}
Let $f \in W^{1,2}(\Tilde{X},T_\lambda)$ be a $\pi_1$-equivariant map,
where $X \in \T(S)$ and $\lambda \in \ML(S)$.  Then
\begin{equation*}
\| \Phi(f) + \frac{1}{4} \phi_X(\lambda) \|_1 \leq \sqrt{2} \left (\E(f) -
\E(\pi_\lambda)\right)^\frac{1}{2} \left ( \E(f)^\frac{1}{2} +
\E(\pi_\lambda)^\frac{1}{2} \right )
\end{equation*}
\end{cor}

\begin{proof}[Proof of Theorem \ref{thm:minimizing}]

Clearly the sequence of harmonic maps $h_i : X \to Y_i$ is a minimizing
sequence.  Applying Theorem \ref{thm:wolf} to the sequence $Y_i$ we
find that the Hopf differentials converge projectively:
$$ \lim_{i\to\infty} [\Phi(h_i)] = [\Phi_\infty], \; \text{where} \; \F(-[\Phi_\infty]) = [\mu].$$

To prove Theorem \ref{thm:minimizing}, we therefore need only show
that the (measurable) Hopf differentials $\Phi(f_i)$ of \emph{any}
minimizing sequence have the same projective limit as the holomorphic
Hopf differentials $\Psi_i = \Phi(h_i)$ of the harmonic maps.

Applying Corollary \ref{cor:energydiff} to such a sequence, we find
\begin{equation}
\label{eqn:ratiolimit}
\begin{split}
 \frac{\| f_i^*(\rho_i) - h_i^*(\rho_i) \|_1}{\E(h_i)} &\leq 
\frac{\sqrt{2} \left ( \E(f_i) - \E(h_i)\right )^\frac{1}{2} 
\left ( \E(f_i)^\frac{1}{2} + \E(h_i)^\frac{1}{2} \right )}{\E(h_i)}\\
&= \sqrt{2} \left ( \frac{\E(f_i)}{\E(h_i)} - 1 \right
)^\frac{1}{2} \left ( 1 + \left ( \frac{\E(f_i)}{\E(h_i)}\right
)^\frac{1}{2} \right ).
\end{split}
\end{equation}
Since $f_i$ is a minimizing sequence, the first factor on the
right hand side of \eqref{eqn:ratiolimit} is $o(1)$ and the second
is bounded, thus
$$ \lim_{i\to\infty} \frac{\| f_i^*(\rho_i) - h_i^*(\rho_i) \|_1}{\E(h_i)} = 0.$$
Finally, the Hopf differential is the $(2,0)$ part of the pullback
metric, so we have
$$ \lim_{i\to\infty} [\Phi(f_i)] = \lim_{i\to\infty} [\Psi_i] = [\Phi], \text{ where }
[\F(-\Phi)] = [\mu].$$
\end{proof}

\section{Proof of the main theorem}

\label{sec:mainproof}

In this section we apply Theorem \ref{thm:minimizing} to the collapsing maps
$\cl_i : X \to Y_i$ to prove Theorem \ref{thm:limit}.

\begin{proof}[Proof of Theorem \ref{thm:limit}]
Fix $X \in \T(S)$ and let $\lambda_i \in \ML(S)$ be a divergent sequence of
measured laminations.  Let $Y_i = \pr_{\lambda_i} X$ so that $X = \gr_{\lambda_i}
Y_i$.  By Lemma \ref{lem:divergence}, the surfaces $Y_i \to \infty$ and thus
$\E(X,Y_i) \to \infty$.  We need to show that if $Y_i \to [\mu] \in \PML(S)$ and
$[\lambda_i] \to [\lambda] \in \PML(S)$ then $i_X([\lambda]) = [\mu]$, or equivalently, that
$[\lambda]$ and $[\mu]$ are the horizontal and vertical measured laminations
of a single holomorphic quadratic differential on $X$.

By Corollary \ref{cor:energyclose},
$$ \E(\cl_i) - \E(X,Y_i) = O(1), $$ while $\E(X,Y_i) \to \infty$, thus
$\cl_i$ is a minimizing sequence (cf. §\ref{sec:convergence}).
Applying Theorem \ref{thm:minimizing}, we conclude that
$$ [\Phi_X(\lambda_i)] \to [\Phi] \in \PQ(X), \text{ where }
[\F(-\Phi)] = [\mu].$$ On the other hand, by Lemma
\ref{lem:intersection}, this implies that $[\F(\Phi)] = [\lambda]$,
and $i_X([\lambda]) = [\mu]$.
\end{proof}

While Theorem \ref{thm:limit} is an asymptotic statement, the ideas
used in the proof above also yield the following finite version of the
comparison between the Hopf differential of the collapsing map and the
holomorphic differential $\phi_X(\lambda)$ representing the grafting
lamination:

\begin{thm}
\label{thm:hopfstrebel}
Let $X \in \T(S)$ and $\lambda \in \ML(S)$. Then the Hopf differential
$\Phi_X(\lambda)$ of the collapsing map $\kappa : X \to \pr_\lambda X$
and the Hubbard-Masur differential $\phi_X(\lambda)$ satisfy
$$ \| 4 \Phi_X(\lambda) - \phi_X(\lambda) \|_1 \leq C\left (1 +
E(\lambda,X)^\frac{1}{2}\right )$$
where $C$ is a constant depending only on $\chi(S)$.
\end{thm}

\begin{proof}
Recall from (\ref{eqn:hopfhopf}) that the Hopf differential of the
co-collapsing map $\hat{\cl}$ is
$$ \Phi(\hat{\cl}) = -\Phi_X(\lambda),$$
while by Lemma \ref{lem:hopfenergyprojection}, that of the harmonic
projection $\pi_\lambda$ is
$$ \Phi(\pi_\lambda) = - \frac{1}{4} \phi_X(\lambda).$$
By Corollary \ref{cor:coenergyclose}, the energy of $\Hat{\cl}$ satisfies
$$ \E(\hat{\cl}) - \E(\pi_\lambda) \leq A := 2 \pi |\chi(S)|.$$
Therefore Corollary \ref{cor:treeenergydiff} implies that
\begin{equation*}
\begin{split}
\| -\Phi_X(\lambda) + \frac{1}{4}\phi_X(\lambda) \|_1 &\leq \sqrt{2}
\left ( \E(\hat{\cl}) - \E(\pi_\lambda) \right)^{\frac{1}{2}} \left (
\E(\pi_\lambda)^{\frac{1}{2}} + \E(\hat{\cl})^{\frac{1}{2}} \right )\\
&\leq \sqrt{2} A^{\frac{1}{2}} 
\left ( \left( \tfrac{1}{2} E(\lambda,X) \right )^{\frac{1}{2}} + 
        \left( \tfrac{1}{2} E(\lambda,X) + A \right
	)^{\frac{1}{2}} \right )\\
& \leq \sqrt{2} A^{\frac{1}{2}} \left ( 2 \left (\tfrac{1}{2}
	E(\lambda,X)\right )^{\frac{1}{2}} + A^{\frac{1}{2}} \right)\\
& = \sqrt{2} A + 2 A^{\frac{1}{2}} E(\lambda,X)^{\frac{1}{2}} 
\end{split}
\end{equation*}
The Theorem then follows by algebra.
\end{proof}

\renewcommand{\thesection}{\Alph{section}}
\setcounter{section}{0}

\section{Appendix: Asymmetry of Teichmüller geodesics}

\label{sec:teichmuller}

Since the antipodal involution $i_X$ is a homeomorphism of $\PML(S)$
to itself, it seems natural to look for an involutive homeomorphism of
$\T(S)$ that has $i_X$ as its boundary values.  In fact, an obvious
candidate is the Teichmüller geodesic involution $I_X : \T(S) \to \T(S)$
that is the pushforward of $(-1) : Q(X) \to Q(X)$ via the Teichmüller
exponential map $\tau : Q(X) \xrightarrow{\sim} \T(S)$.  However, we now
sketch an example showing that $I_X$ does not extend continuously to
the Thurston compactification of $\T(S)$, leading us to view the
pruning map based at $X$, $\lambda \mapsto \pr_\lambda X$, as a kind of
substitute for the Teichmüller geodesic involution that \emph{does}
extend continuously to the antipodal map of $\PML$.

\begin{figure}
\begin{center}
% -- BEGIN "lratios.inc" --
\begin{picture}(0,0)%
\includegraphics{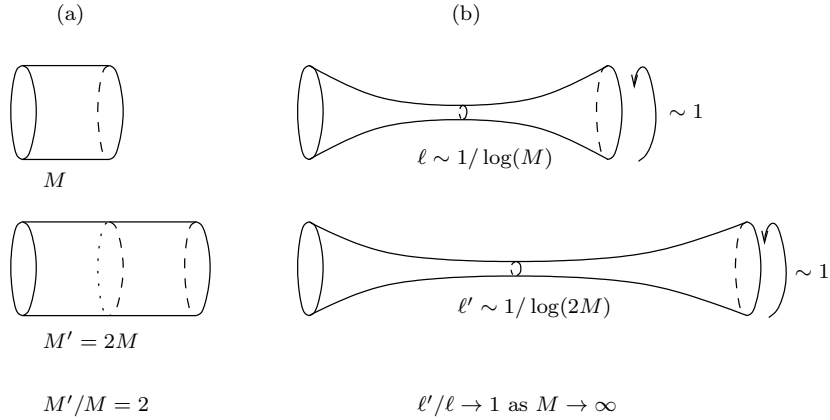}%
\end{picture}%
\setlength{\unitlength}{4144sp}%
\begingroup\makeatletter\ifx\SetFigFont\undefined%
\gdef\SetFigFont#1#2#3#4#5{%
  \reset@font\fontsize{#1}{#2pt}%
  \fontfamily{#3}\fontseries{#4}\fontshape{#5}%
  \selectfont}%
\fi\endgroup%
\begin{picture}(4718,2529)(2343,-2407)
\put(7046,-1552){\makebox(0,0)[lb]{\smash{{\SetFigFont{8}{9.6}{\familydefault}{\mddefault}{\updefault}{$\sim 1$}%
}}}}
\put(6292,-590){\makebox(0,0)[lb]{\smash{{\SetFigFont{8}{9.6}{\familydefault}{\mddefault}{\updefault}{$\sim 1$}%
}}}}
\put(2627, -3){\makebox(0,0)[lb]{\smash{{\SetFigFont{8}{9.6}{\familydefault}{\mddefault}{\updefault}{(a)}%
}}}}
\put(4992, -3){\makebox(0,0)[lb]{\smash{{\SetFigFont{8}{9.6}{\familydefault}{\mddefault}{\updefault}{(b)}%
}}}}
\put(5030,-1762){\makebox(0,0)[lb]{\smash{{\SetFigFont{8}{9.6}{\familydefault}{\mddefault}{\updefault}{$\ell' \sim 1/\log(2 M)$}%
}}}}
\put(4789,-865){\makebox(0,0)[lb]{\smash{{\SetFigFont{8}{9.6}{\familydefault}{\mddefault}{\updefault}{$\ell \sim 1/\log(M)$}%
}}}}
\put(4789,-2348){\makebox(0,0)[lb]{\smash{{\SetFigFont{8}{9.6}{\familydefault}{\mddefault}{\updefault}{$\ell' / \ell \rightarrow 1$ as $M \rightarrow \infty$}%
}}}}
\put(2548,-2348){\makebox(0,0)[lb]{\smash{{\SetFigFont{8}{9.6}{\familydefault}{\mddefault}{\updefault}{$M'/M = 2$}%
}}}}
\put(2548,-1968){\makebox(0,0)[lb]{\smash{{\SetFigFont{8}{9.6}{\familydefault}{\mddefault}{\updefault}{$M' = 2M$}%
}}}}
\put(2548,-1003){\makebox(0,0)[lb]{\smash{{\SetFigFont{8}{9.6}{\familydefault}{\mddefault}{\updefault}{$M$}%
}}}}
\end{picture}%
% -- END "lratios.inc" --
\caption{(a) Euclidean cylinders with moduli $M$ and $2M$ correspond
 to (b) hyperbolic cylinders whose core geodesics have approximately the
 same length.  This phenomenon leads to Teichmüller rays for
 distinct Strebel differentials that converge to the same point in the
 Thurston boundary of Teichmüller space \cite{Masur}.}
\label{fig:lengths}
\end{center}
\end{figure}

By a theorem of Masur, if
$\phi$ is a Strebel differential on $X$ whose trajectories
represent homotopy classes $(\alpha_1, \ldots ,\alpha_n)$, then the Teichmüller ray
determined by $\phi$ converges to the point $[\alpha_1 + \cdots + \alpha_n] \in \PML(S)$ in
the Thurston compactification \cite{Masur}.  Note that the limit point corresponds
to a measured lamination in which each curve $\alpha_i$ has the same
weight, independent of the relative sizes of the cylinders on $X$
determined by $\phi$.  This happens because the Thurston boundary
reflects the geometry of hyperbolic geodesics on the surface, and
hyperbolic length is approximated by the reciprocal of the logarithm
of a cylinder's height, as in Figure \ref{fig:lengths}.

Now suppose $\phi$ and $\psi$ are holomorphic quadratic differentials on a
Riemann surface $X$ such that each of $± \phi, ± \psi$ is Strebel,
where the trajectories of $\phi$ and $\psi$ represent different sets of
homotopy classes, while those of $-\phi$ and $-\psi$ represent the same
homotopy classes.  Then by Masur's theorem, the Teichmüller geodesics
corresponding to $\phi$ and $\psi$ converge to the same point on $\PML(S)$ in
the negative direction, while in the positive direction they converge
to distinct points.  If $I_X$ were to extend to a continuous map of
the Thurston boundary, then any pair of Teichmüller geodesics that are
asymptotic in one direction would necessarily be asymptotic in both
directions, thus no such extension exists.

\begin{figure}
\begin{center}
\includegraphics[width=8cm]{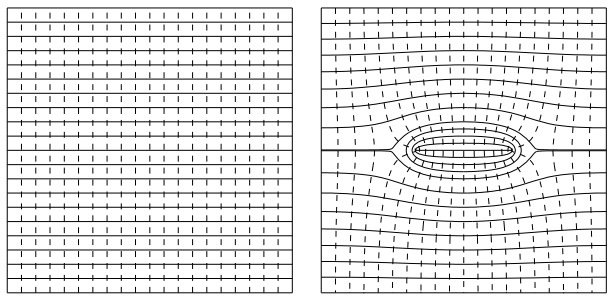}
\caption{Horizontal (solid) and vertical (dashed) trajectories of
  $\phi_0$ and $\psi_0$ on the square torus.}
\label{fig:trajectories}
\end{center}
\end{figure}

One can explicitly construct such $(X,\phi,\psi)$ as follows: Let $X_0$
denote the square torus $\C / (2\Z \oplus 2 i \Z)$, and let $\psi_0 = \psi_0(z)
dz^2$ be a meromorphic quadratic differential on $X_0$ with simple
zeros at $z = ± \frac{1 + \epsilon}{2}$ and simple poles at $z = ±
\frac{1}{2}$ and such that $\psi_0(x) \in \R$ for $x \in \R$.  Such a
differential exists by the Abel-Jacobi theorem, and in fact is given
by
$$ \psi_0(z) dz^2 = \frac{\wp(z) - c_0}{\wp(z) - c_1}$$ for suitable
constants $c_i$, where $\wp(z)$ is the Weierstrass function for $X_0$.
Let $\phi_0 = dz^2$, a holomorphic quadratic differential on $X_0$.
The trajectory structures of $\phi_0$ and $\psi_0$ are displayed in
Figure \ref{fig:trajectories}.

Let $X$ be the surface of genus $2$ obtained as a 2-fold cover of
$X_0$ branched over $± \frac{1}{2}$; then $\psi_0$ and $\phi_0$ determine
holomorphic quadratic differentials $\psi$ and $\phi$ on $X$, where the lift
of $\psi_0$ is holomorphic because the simple poles at $± \frac{1}{2}$
are branch points of the covering map $X \to X_0$.

\begin{figure}
\begin{center}
% -- BEGIN "torusglue.inc" --
\begin{picture}(0,0)%
\includegraphics{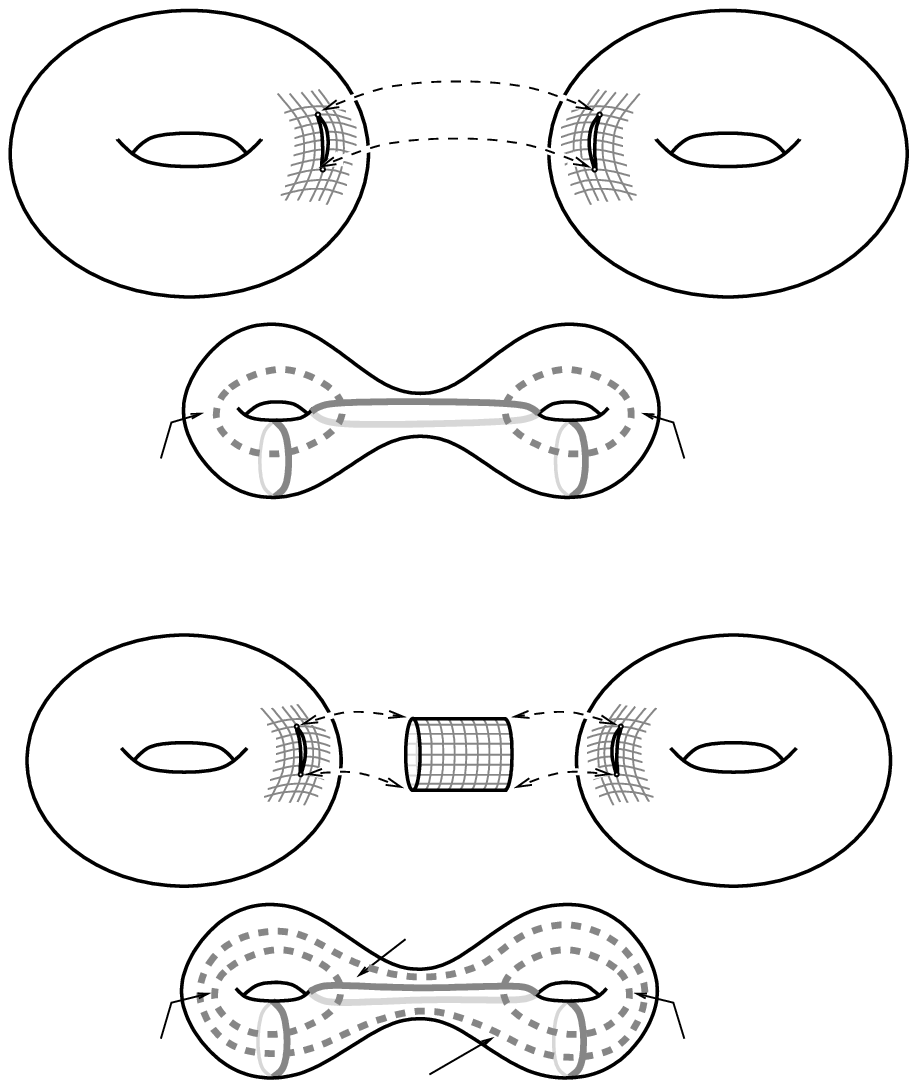}%
\end{picture}%
\setlength{\unitlength}{4144sp}%
\begingroup\makeatletter\ifx\SetFigFont\undefined%
\gdef\SetFigFont#1#2#3#4#5{%
  \reset@font\fontsize{#1}{#2pt}%
  \fontfamily{#3}\fontseries{#4}\fontshape{#5}%
  \selectfont}%
\fi\endgroup%
\begin{picture}(4524,5147)(2306,-6425)
\put(4529,-6249){\makebox(0,0)[lb]{\smash{{\SetFigFont{8}{9.6}{\familydefault}{\mddefault}{\updefault}{$\alpha$}%
}}}}
\put(2321,-1998){\makebox(0,0)[lb]{\smash{{\SetFigFont{9}{10.8}{\familydefault}{\mddefault}{\updefault}{a)}%
}}}}
\put(2321,-4782){\makebox(0,0)[lb]{\smash{{\SetFigFont{9}{10.8}{\familydefault}{\mddefault}{\updefault}{b)}%
}}}}
\put(5198,-6367){\makebox(0,0)[lb]{\smash{{\SetFigFont{8}{9.6}{\familydefault}{\mddefault}{\updefault}{$\eta_-$}%
}}}}
\put(3832,-6367){\makebox(0,0)[lb]{\smash{{\SetFigFont{8}{9.6}{\familydefault}{\mddefault}{\updefault}{$\eta_+$}%
}}}}
\put(4556,-5537){\makebox(0,0)[lb]{\smash{{\SetFigFont{8}{9.6}{\familydefault}{\mddefault}{\updefault}{$\beta$}%
}}}}
\put(3329,-6145){\makebox(0,0)[lb]{\smash{{\SetFigFont{8}{9.6}{\familydefault}{\mddefault}{\updefault}{$\gamma_+$}%
}}}}
\put(5793,-6145){\makebox(0,0)[lb]{\smash{{\SetFigFont{8}{9.6}{\familydefault}{\mddefault}{\updefault}{$\gamma_-$}%
}}}}
\put(3862,-3697){\makebox(0,0)[lb]{\smash{{\SetFigFont{8}{9.6}{\familydefault}{\mddefault}{\updefault}{$\eta_+$}%
}}}}
\put(5220,-3697){\makebox(0,0)[lb]{\smash{{\SetFigFont{8}{9.6}{\familydefault}{\mddefault}{\updefault}{$\eta_-$}%
}}}}
\put(4559,-2874){\makebox(0,0)[lb]{\smash{{\SetFigFont{8}{9.6}{\familydefault}{\mddefault}{\updefault}{$\beta$}%
}}}}
\put(3314,-3499){\makebox(0,0)[lb]{\smash{{\SetFigFont{8}{9.6}{\familydefault}{\mddefault}{\updefault}{$\gamma_+$}%
}}}}
\put(5756,-3492){\makebox(0,0)[lb]{\smash{{\SetFigFont{8}{9.6}{\familydefault}{\mddefault}{\updefault}{$\gamma_-$}%
}}}}
\end{picture}%
% -- END "torusglue.inc" --
\caption{Two constructions of a surface of genus two with a Strebel
  differential:  (a) Two tori are glued along a segment of a leaf of a
  Euclidean foliation. (b) Two tori are glued to the ends of a
  Euclidean cylinder.  Below each example, the homotopy classes
  represented by the horizontal and vertical trajectories are shown
  (as solid and dashed lines, respectively).  When the tori and
  cylinder are chosen correctly, this construction produces an example
  of Teichmüller geodesics that are asymptotic in one direction
  while converging to distinct endpoints in the opposite direction.}
\label{fig:strebel_curves}
\end{center}
\end{figure}

Both $\phi$ and $\psi$ have closed vertical and horizontal trajectories as
in the construction above (see Figure \ref{fig:strebel_curves}).
Specifically, let $\gamma$ and $\eta$ denote the free homotopy classes of
simple closed curves on $X_0$ that arise as the quotients of $\R$ and
$i \R$, respectively; both $\gamma$ and $\eta$ have two distinct lifts
($\gamma_{±}$ and $\eta_{±}$, respectively) to $X$.  Let $\alpha$ denote the
separating curve on $X$ that covers $[-1/2,1/2]$, and let $\beta$ denote
the simple closed curve on $X$ that is the union of the two lifts of
$[-i,i]$.  Then:
\begin{enumerate}
\item the trajectories of $\phi$ represent $(\gamma_+,\gamma_-)$,
\item the trajectories of $\psi$ represent $(\gamma_+,\gamma_-,\alpha)$, and
\item the trajectories of both $-\phi$ and $-\psi$ represent $(\eta_+,\eta_-,\beta)$.
\end{enumerate}

\nocite{}
%% -- BEGIN BIBLIOGRAPHY --

%% -- END BIBLIOGRAPHY --


\begin{thebibliography}{Mas2}

\bibitem[Bes]{Be}
Mladen Bestvina.
\newblock {Degenerations of the hyperbolic space}.
\newblock {\em Duke Math. J.} {\bf 56}(1988), 143--161.

\bibitem[CS]{CS}
Marc Culler and Peter~B. Shalen.
\newblock {Varieties of group representations and splittings of
  {$3$}-manifolds}.
\newblock {\em Ann. of Math. (2)} {\bf 117}(1983), 109--146.

\bibitem[D]{Du}
David Dumas.
\newblock {\em Complex Projective Structures, Grafting, and Teichm\"uller
  Theory}.
\newblock PhD thesis, Harvard University, May 2004.

\bibitem[FW]{FW}
Benson Farb and Michael Wolf.
\newblock {Harmonic splittings of surfaces}.
\newblock {\em Topology} {\bf 40}(2001), 1395--1414.

\bibitem[FLP]{FLP}
A.~Fathi, F.~Laudenbach, and V.~Poenaru.
\newblock {\em Travaux de Thurston sur les surfaces}, volume~66 of {\em
  Ast\'erisque}.
\newblock Soci\'et\'e Math\'ematique de France, Paris, 1979.
\newblock S\'eminaire Orsay, With an English summary.

\bibitem[Gar]{Ga}
Frederick~P. Gardiner.
\newblock {\em Teichm\"uller theory and quadratic differentials}.
\newblock Pure and Applied Mathematics. John Wiley \& Sons Inc., New York,
  1987.

\bibitem[Har]{Ha}
Philip Hartman.
\newblock {On homotopic harmonic maps}.
\newblock {\em Canad. J. Math.} {\bf 19}(1967), 673--687.

\bibitem[HM]{HM}
John Hubbard and Howard Masur.
\newblock {Quadratic differentials and foliations}.
\newblock {\em Acta Math.} {\bf 142}(1979), 221--274.

\bibitem[KT]{KT}
Yoshinobu Kamishima and Ser~P. Tan.
\newblock {Deformation spaces on geometric structures}.
\newblock In {\em Aspects of low-dimensional manifolds}, pages 263--299.
  Kinokuniya, Tokyo, 1992.

\bibitem[Kap]{Ka}
Michael Kapovich.
\newblock {\em Hyperbolic manifolds and discrete groups}.
\newblock Birkh\"auser Boston Inc., Boston, MA, 2001.

\bibitem[KS1]{KS}
Nicholas~J. Korevaar and Richard~M. Schoen.
\newblock {Sobolev spaces and harmonic maps for metric space targets}.
\newblock {\em Comm. Anal. Geom.} {\bf 1}(1993), 561--659.

\bibitem[KS2]{KS2}
Nicholas~J. Korevaar and Richard~M. Schoen.
\newblock {Global existence theorems for harmonic maps to non-locally compact
  spaces}.
\newblock {\em Comm. Anal. Geom.} {\bf 5}(1997), 333--387.

\bibitem[KP]{KP2}
Ravi~S. Kulkarni and Ulrich Pinkall.
\newblock {A canonical metric for {M}\"obius structures and its applications}.
\newblock {\em Math. Z.} {\bf 216}(1994), 89--129.

\bibitem[Lev]{Lev}
Gilbert Levitt.
\newblock {Foliations and laminations on hyperbolic surfaces}.
\newblock {\em Topology} {\bf 22}(1983), 119--135.

\bibitem[Mas1]{Masur2}
Howard Masur.
\newblock {On a class of geodesics in {T}eichm\"uller space}.
\newblock {\em Ann. of Math. (2)} {\bf 102}(1975), 205--221.

\bibitem[Mas2]{Masur}
Howard Masur.
\newblock {Two boundaries of {T}eichm\"uller space}.
\newblock {\em Duke Math. J.} {\bf 49}(1982), 183--190.

\bibitem[Min]{Mi}
Yair~N. Minsky.
\newblock {Harmonic maps, length, and energy in {T}eichm\"uller space}.
\newblock {\em J. Differential Geom.} {\bf 35}(1992), 151--217.

\bibitem[MS]{MS}
John~W. Morgan and Peter~B. Shalen.
\newblock {Valuations, trees, and degenerations of hyperbolic structures. {I}}.
\newblock {\em Ann. of Math. (2)} {\bf 120}(1984), 401--476.

\bibitem[Ota]{Otal}
Jean-Pierre Otal.
\newblock {Le th\'eor\`eme d'hyperbolisation pour les vari\'et\'es fibr\'ees de
  dimension 3}.
\newblock {\em Ast\'erisque} {\bf 235}(1996).

\bibitem[Pau]{Pa}
Fr{\'e}d{\'e}ric Paulin.
\newblock {Topologie de {G}romov \'equivariante, structures hyperboliques et
  arbres r\'eels}.
\newblock {\em Invent. Math.} {\bf 94}(1988), 53--80.

\bibitem[SW]{SW}
Kevin~P. Scannell and Michael Wolf.
\newblock {The grafting map of {T}eichm\"uller space}.
\newblock {\em J. Amer. Math. Soc.} {\bf 15}(2002), 893--927.
\newblock  {{\tt arXiv:math.DG/9810082}}

\bibitem[Sko]{Skora}
Richard~K. Skora.
\newblock {Splittings of surfaces}.
\newblock {\em J. Amer. Math. Soc.} {\bf 9}(1996), 605--616.

\bibitem[Tan]{Ta}
Harumi Tanigawa.
\newblock {Grafting, harmonic maps and projective structures on surfaces}.
\newblock {\em J. Differential Geom.} {\bf 47}(1997), 399--419.
\newblock  {{\tt arXiv:math.DG/9508216}}

\bibitem[Thu]{Th2}
William~P. Thurston.
\newblock {Geometry and Topology of Three-Manifolds}.
\newblock Princeton lecture notes, 1979.
\newblock  {{\tt http://www.msri.org/publications/books/gt3m/}}

\bibitem[W1]{W}
Michael Wolf.
\newblock {The {T}eichm\"uller theory of harmonic maps}.
\newblock {\em J. Differential Geom.} {\bf 29}(1989), 449--479.

\bibitem[W2]{W4}
Michael Wolf.
\newblock {Harmonic maps from a surface and degeneration in {T}eichm\"uller
  space}.
\newblock In {\em Low-dimensional topology (Knoxville, TN, 1992)}, Conf. Proc.
  Lecture Notes Geom. Topology, III, pages 217--239. Internat. Press,
  Cambridge, MA, 1994.

\bibitem[W3]{W5}
Michael Wolf.
\newblock {Harmonic maps from surfaces to {$\bold R$}-trees}.
\newblock {\em Math. Z.} {\bf 218}(1995), 577--593.

\bibitem[W4]{W6}
Michael Wolf.
\newblock {On realizing measured foliations via quadratic differentials of
  harmonic maps to {$\bold R$}-trees}.
\newblock {\em J. Anal. Math.} {\bf 68}(1996), 107--120.

\bibitem[W5]{W7}
Michael Wolf.
\newblock {Measured foliations and harmonic maps of surfaces}.
\newblock {\em J. Differential Geom.} {\bf 49}(1998), 437--467.

\end{thebibliography}
\end{document}